\documentclass[11pt]{article}

\usepackage[active]{srcltx}
\usepackage[OT2,T1]{fontenc}
\usepackage{bbm}\usepackage{amsfonts}\usepackage{amssymb}\usepackage{amsmath,amsthm}

\DeclareSymbolFont{cyrletters}{OT2}{wncyr}{m}{n}\DeclareSymbolFont{cyrlettersb}{OT2}{wncyr}{b}{n}
\DeclareMathSymbol{\Sha}{\mathalpha}{cyrletters}{"58}\DeclareMathSymbol{\Ya}{\mathalpha}{cyrlettersb}{"17}

\newcommand{\A}[3]{\big\{\!\! \begin{smallmatrix} #1 & #2\\ #3 & \end{smallmatrix}\!\!\big\} }
\newcommand{\Au}[4]{\big\{\!\! \begin{smallmatrix} #1 & #2\\ #3 & \end{smallmatrix}\!\!\big\} ^{\uparrow #4}}
\newcommand{\x}[1]{x_{#1\rangle}^{\phantom{\uparrow}}}
\newcommand{\y}[1]{y_{#1\rangle}^{\phantom{\uparrow}}}
\newcommand{\z}[1]{z_{#1\rangle}^{\phantom{\uparrow}}}
\newcommand{\xu}[2]{x_{#1\rangle}^{\uparrow #2}}\newcommand{\yu}[2]{y_{#1\rangle}^{\uparrow #2}}
\newcommand{\zu}[2]{z_{#1\rangle}^{\uparrow #2}}
\renewcommand{\k}{\mathfrak{k}}\newcommand{\id}{\ensuremath{\mathbbm{1}}}

\def\be{\begin{equation}}\def\ee{\end{equation}}\def\ba{\begin{eqnarray}}\def\ea{\end{eqnarray}}
\def\lb{\label}
\def\s{\sigma}\def\b{\beta}\def\o{\omega}

\begin{document}
\begin{center}

\hfill \small{CPT-P03-2008}

\vskip .2cm

{\Large\bf Braidings of Tensor Spaces}\end{center}
\vspace{.3cm}

\begin{center}{\large\bf T. Grapperon\footnote{grapperon@cpt.univ-mrs.fr} and O.V. Ogievetsky\footnote{oleg@cpt.univ-mrs.fr}}\end{center}

\begin{center}
Center of Theoretical Physics\footnote{Unit\'e Mixte de Recherche
(UMR 6207) du CNRS et des Universit\'es Aix--Marseille I,
Aix--Marseille II et du Sud Toulon -- Var; laboratoire affili\'e \`a la
FRUMAM (FR 2291)}, Luminy,
13288 Marseille, France
\end{center}

\vspace{.6 cm}
\centerline{\footnotesize\bf Abstract} 
{\footnotesize\noindent 
\begin{center}
\begin{minipage}[c]{11.4cm}
Let $V$ be a braided vector space, that is, a vector space together with a solution 
$\hat{R}\in {\text{End}}(V\otimes V)$ of the Yang--Baxter equation. Denote
$T(V):=\bigoplus_k V^{\otimes k}$. We associate to $\hat{R}$ a solution 
$T(\hat{R})\in {\text{End}}(T(V)\otimes T(V))$ of the Yang--Baxter equation
on the tensor space $T(V)$. The correspondence $\hat{R}\rightsquigarrow T(\hat{R})$ 
is functorial with respect to $V$. 
\end{minipage}\end{center}}

\vspace{.2cm}
\section{Introduction}

Solutions of the Yang--Baxter equation find applications in statistical models, integrable systems, knot theory, representation
theory of braid groups and many other areas. In the present article we propose a method for constructing new solutions of 
the Yang--Baxter equation. These new solutions live in the tensor space $T(V):=\bigoplus_k V^{\otimes k}$, where $V$ is
a vector space equipped with a solution $\hat{R}$ of the Yang--Baxter equation (in other words, $V$ is a {\it braided} vector
space). The initial data for the new solution includes, in addition to the operator $\hat{R}$, a scalar parameter $q$.  

We make a "block triangular" (in the sense specified in the section \ref{braidingtensors}) Ansatz. It turns out that there is a "universal" 
system of equations in the algebra of the infinite Artin braid group (see section \ref{notation-convention} for the definition of 
the braid group and connections between braid groups and braided vector spaces), solving which we automatically produce the solution 
of the Yang--Baxter equation on $T(V)$; we present a solution of this universal system. As a by-product of the universality, the 
correspondence $\hat{R}\rightsquigarrow T(\hat{R})$ is functorial in the following sense. Let $V$ and $V'$ be two braided 
vector spaces, with braidings $\hat{R}$ and $\hat{R}'$ respectively. Let $f\in \text{Hom}(V,V')$ be a homomorhism of vector spaces; 
it extends to the homomorphism $T(f): T(V)\rightarrow T(V')$ by the natural rule: the restriction of $T(f)$ on $V^{\otimes k}$ is 
$f\otimes f\otimes\dots\otimes f$ ($k$ times). The functoriality means: if $f$ intertwines $\hat{R}$ and $\hat{R}'$, that is, 
$(f\otimes f) \hat{R}=\hat{R}' (f\otimes f)$, then $T(f)$ intertwines  $T(\hat{R})$ and $T(\hat{R}')$.

Given a braiding on a vector spave $V$, there is a well-known standard braiding on $T(V)$ which we call diagonal (see sections 
\ref{notation-convention} and \ref{braidingtensors} for details). It should be noted that the braiding $T(\hat{R})$ is different from the 
standard diagonal braiding: it mixes tensors with different number of indices.

The solution of the universal system is formulated with the help of two kinds of elements in the braid group algebra: generalizations 
of the binomial coefficients, called {\it braid shuffle} elements, and the {\it braid Pochhammer symbols}; see section \ref{notation-convention}
for definitions. The braid shuffle elements and Pochhammer symbols satisfy certain identities, generalizing the corresponding classical identities.
These generalized identities (some of them appear to be new and are given with a proof) are presented in section \ref{notation-convention}
as well. The solution itself is given in section \ref{braidingtensors}; also, some additional 
properties of the solution $T(\hat{R})$ and some directions of future research are outlined there.

\section{Notation, conventions and preliminary identities}
\label{notation-convention}

Formulas hereafter often contain sums and products. It is always assumed that an empty sum (that is, $\sum_{j=a}^{a-1}$ for any integer $a$) 
is equal to zero and an empty product (that is, $\prod_{j=a}^{a-1}$ for any integer $a$) is equal to 1. 
The set of non-negative integers is denoted by ${\mathbb Z}_{\geq 0}$. 

\vskip .2cm
\noindent {\bf Braid group.} We shall use the Artin presentation of the braid group $B_{n+1}$ by generators $\sigma_i$, $1\leq i\leq n$, and relations 
\begin{eqnarray}\label{ryb1}\sigma_{i}\sigma_{j}\sigma_{i}=\sigma_{j}\sigma_{i}\sigma_{j}
&{\mathrm{if}}& |i-j|=1\ ,\\[.5em]
\sigma_{i}\sigma_{j}=\sigma_{j}\sigma_{i}&{\mathrm{if}}&|i-j|>1\ .\label{ryb2}\end{eqnarray}
The first of these is called the \emph{braid} relation and the second one the \emph{far 
commutativity} (it reflects the commutativity of unlinked braids). 

One has the {\it tower} $B_0\subset B_1\subset B_2\subset\cdots\subset B_\infty$ (the groups $B_0$ and $B_1$ are trivial) 
defined by inclusions $B_n\ni\sigma_i\mapsto \sigma_i\in B_{n+1}$, $i=1,\dots ,n-1$;
the group $B_\infty =\displaystyle{\lim_{\longrightarrow}}\, B_n$ is the inductive limit,
the group with the countable set $\{\sigma_1 ,\sigma_2 ,\dots\}$ of generators subject to the above relations.

Given $\ell\in {\mathbb Z}_{\geq 0}$, we denote the endomorphism of $B_\infty$, sending $\sigma_i$ to $\sigma_{i+\ell}$, $i=1,2,\dots$, by 
$\ ^{\uparrow\ell}\,$, 
\begin{equation} ^{\uparrow\ell}:\sigma_i\mapsto\sigma_{i+\ell} \ ,\end{equation} 
as in \cite{OP}. We write $w^{\uparrow\ell}$ for the image of an element $w$; for example, $\sigma_i^{\uparrow\ell}=\sigma_{i+\ell}$.

\vskip .2cm
\noindent {\bf Braided vector spaces.} 
Let $V$ a braided vector space of dimension $N$ over a field $\k$, with a braiding $\hat{R}\in\mathrm{Aut}(V^{\otimes 2})$; $\hat{R}$ is 
an invertible solution of the Yang--Baxter equation. 
The braiding induces the so-called {\it local} representation $\rho_{\hat{R}}$ of the braid 
group tower, $\rho_{\hat{R}} :\, B_n\rightarrow \mathrm{Aut}(V^{\otimes n})$, by $\rho_{\hat{R}} (\sigma_i)=\hat{R}_i:= 
\id_V^{\otimes (i-1)}\otimes \hat{R}\otimes\id_V^{\otimes(n-i-1)}$ (that is, 
$\sigma_i$ acts as $\hat{R}$ on the copies number $i$ and $i+1$ of the space $V$ in $V\otimes V\otimes V...$ and 
as the identity on the other copies). We use the same notation $\ ^{\uparrow j}\,$ for the shift
in the copies of the space $V$; for example, $\rho_{\hat{R}} (\sigma_i)=\rho_{\hat{R}} 
(\sigma_1)^{\uparrow i-1}$.

In the sequel we slightly abuse the notation writing $\sigma$ instead of $\rho_{\hat{R}}  (\sigma)$; this should not produce any 
confusion: one can say that we consider $V\otimes V\otimes V...$ as a $B_\infty$-module, $\sigma_i$ acts as $\hat{R}_i$.
When needed we shall specify the Yang-Baxter operator $\hat{R}$ defining the local representation.

We use the notation $\x{a}$, as in \cite{IO}, for a vector in the copy number $a$ of the space
$V$; $[a]$ is the set $\{ 1,2,\dots ,a\}$; we write $\x{[a]}$ for a tensor 
$x^{i_1 i_2 \dots i_a}\in V^{\otimes a}$ with $a$ indices and 
$\xu{[a]}{b}$ for the same tensor $x$ belonging to the copies numbered from $1+b$ to $a+b$;
for example, $\x{[a]}\yu{[b]}{a}$ is an element $x^{i_1 i_2 \dots i_a}\otimes y^{j_1 j_2 \dots j_b}\in V^{\otimes a}\otimes V^{\otimes b}$.

The braiding can be understood as a rule of exchanging two copies of the space $V$,
\be\lb{br1} x_{1\rangle}y_{2\rangle}=\sigma_1\ y_{1\rangle}x_{2\rangle}\ .\ee
Assume the same exchange rule for copies $i$ and $i+1$ ($i$ is arbitrary) of the space $V$,
$x_{i\rangle}y_{i+1\rangle}=\sigma_i\ y_{i\rangle}x_{i+1\rangle}$; assume the same exchange rules for 
$x_{i\rangle}z_{i+1\rangle}$ and $y_{i\rangle}z_{i+1\rangle}$. The Yang-Baxter relation for $\hat{R}$ 
(the braid relation for $\sigma$) ensures the coincidence of rewriting $x_{1\rangle}y_{2\rangle}z_{3\rangle}$ in two ways 
(starting with $x_{1\rangle}y_{2\rangle}$ or $y_{2\rangle}z_{3\rangle}$) to the form 
$z_{1\rangle}y_{2\rangle}x_{3\rangle}$. 

The arrangement of the product $\x{1}\dots\x{k}\y{k+1}\dots\y{k+\ell}$ to the form 
$\y{1}\dots\y{\ell}\x{\ell +1}\dots\x{\ell +k}$ induces the exchange rule for 
decomposable tensors,
\be\lb{br2}\x{1}\x{2}\dots\x{k}\ \y{k+1}\y{k+2}\dots\y{k+\ell}=\b_{k,l}\ 
\y{1}\y{2}\dots\y{\ell}\ \x{\ell +1}\x{\ell +2}\dots\x{\ell +k}\ ,\ee
$\beta_{k,l}: V^{\otimes k}\otimes V^{\otimes l}\rightarrow V^{\otimes l}\otimes V^{\otimes k}$.
For example, $\b_{k,1}=\sigma_k\sigma_{k-1}\dots\sigma_1$ and $\b_{1,l}=\sigma_1\sigma_2\dots\sigma_l$. 
The arrangement can be done, among other ways, by moving $y$'s one after another to the left 
or $x$'s to the right, giving two expessions for $\b_{k,l}$,
\be\label{braiding-diag}\begin{array}{rcl}\beta_{k,l}&=&(\sigma_k\sigma_{k+1}\dots\sigma_{k+l-1})
(\sigma_{k-1}\sigma_k\dots\sigma_{k+l-2})\dots (\sigma_1\sigma_2\dots\sigma_l)\\[.5em]
&=&(\sigma_k\sigma_{k-1}\dots\sigma_1)(\sigma_{k+1}\sigma_k\dots\sigma_2)\dots
(\sigma_{k+l-1}\sigma_{k+l-2}\dots\sigma_l)\ .\end{array}\ee
The equality of two expressions for $\beta_{kl}$ in (\ref{braiding-diag}) is guaranteed by the Yang--Baxter equation 
for $\sigma$. 

The identities 
\be\lb{br4}\b_{l+m,n}^{\phantom{\uparrow}}=\b_{m,n}^{\uparrow l}\b_{l,n}^{\phantom{\uparrow}}
\quad ,\quad \b_{l,m+n}^{\phantom{\uparrow}}=
\b_{l,m}^{\phantom{\uparrow}}\b_{l,n}^{\uparrow m}\ee
ensure the consistency of eq.(\ref{br2}) and it then follows that 
\be\lb{br5}\b_{m,j}^{\phantom{\uparrow}}\b_{m+l,n}^{\uparrow j}=\b_{l,n}^{\uparrow m+j}
\b_{m,j+n}^{\phantom{\uparrow}}\quad ,\quad \b_{n,m+l}^{\uparrow j}\b_{j,m}^{\phantom{\uparrow}}
=\b_{j+n,m}^{\phantom{\uparrow}}\b_{n,l}^{\uparrow m+j}\ .\ee
We shall often use the identity
\be\lb{br3}\b_{k,l}\phi\psi^{\uparrow l}=\psi\phi^{\uparrow k}\b_{k,l}\ ,\ee
valid for arbitrary elements $\phi\in\k B_l$ and $\psi\in\k B_k$ of the group algebras.

\vskip .2cm
\noindent {\bf Lift of the longest element.} Denote by $\o_a$ the lift of the longest element from the symmetric group 
$S_a$ to $B_a$. It satisfies
\be\lb{br6}\o_a\phi =\phi'\o_a\ \ \ \forall\ \phi\in\k B_a\ ,\ee
where $':B_a\to B_a$ is the automorphism $\s_i\mapsto\s_{a-i}$ of $B_a$. We have
\be\lb{br7}\o_{a+b}^{\phantom{\uparrow}}=\b_{a,b}^{\phantom{\uparrow}}\ \o_a^{\uparrow b}\o_b^{\phantom{\uparrow}}\ .\ee

\vskip .2cm
\noindent {\bf Shuffles.} We shall use braid analogues of the ($q$-)binomial coefficients, the braid shuffle elements $\Sha_{m,n}$, 
indexed by two integers $m,n$ with $m+n\in {\mathbb Z}_{\geq 0}$. The shuffle element $\Sha_{m,n}$ 
belongs to the group ring $\k B_{m+n}$ (and therefore to the inductive limit of the tower, $\k B_\infty$); the
shuffle elements are uniquely defined by the initial condition 
$\Sha_{n,-n}=\delta_{n,0}$, $n\in  {\mathbb Z}$, ($\delta$ is the Kronecker symbol)
and any of the recurrence braid analogues of the Pascal rule (these two recurrences produce the same sets of
shuffle elements)
\begin{eqnarray}\lb{shufrec}\Sha_{m,n}^{\phantom{\uparrow}}&=&\Sha_{m-1,n}^{\phantom{\uparrow}}+
\Sha_{m,n-1}^{\phantom{\uparrow}}\beta_{m,1}^{\uparrow n-1}\ ,\\[.5em]
\lb{shufrec2}\Sha_{m,n}^{\phantom{\uparrow}}&=&
\Sha_{m,n-1}^{\uparrow 1}+\Sha_{m-1,n}^{\uparrow 1}\beta_{1,n}^{\phantom{\uparrow}}\ .\end{eqnarray} 
In particular, $\Sha_{0,n}=1$ and $\Sha_{n,0}=1$ for 
$n\in {\mathbb Z}_{\geq 0}$; for $i+j\in {\mathbb Z}_{\geq 0}$, $\Sha_{i,j}=0$ if
$i<0$ or $j<0$. We have 
\begin{equation}\Sha_{n+k,m}\Sha_{k,n}^{\uparrow m}=\Sha_{k,m+n}\Sha_{n,m}\ .\label{shaiden}\end{equation}
The shuffle elements are denoted in several different ways in the literature. We use the notation of \cite{IO2}. 

\vskip .2cm
\noindent {\bf Braid Pochhammer symbols.} Let 
\begin{equation} P_{k,n}(x,y):=(x-\beta_{k,1}y)(x-\beta_{k+1,1}y)\dots (x-\beta_{k+n-1,1}y)\ .
\end{equation}
Here $x$ and $y$ are parameters. The elements $P_{k,n}(x,y)$ (defined for $k,n\in {\mathbb Z}_{\geq 0}$) can be understood as braid analogues of the  ($q$-)Pochhammer symbols.

By definition,
\begin{equation} P_{k,n}(x,y)=P_{k,a}(x,y)P_{k+a,n-a}(x,y)\ .\label{recpo}\end{equation}

\noindent {\bf Lemma.} We have
\begin{equation} P_{k,n}(x,z)=\sum_a\Sha_{n-a,a}^{\uparrow k}\beta_{k,a}P_{0,a}(y,z)P_{k,n-a}(x,y)^{\uparrow a}\ .
\label{robrbin}\end{equation}

\noindent {\tt Note.} In (\ref{robrbin}) and in sums hereafter we often omit the range of summation: one may understand 
the range as $ {\mathbb Z}$; however the sums are finite because the summands differ from 0 only for a finite number of
values; for instance, in (\ref{robrbin}) the shuffle element differs from zero for $a$ from 0 to $n$ only. 

\vskip .2cm
\noindent {\it Proof.} Induction in $n$. For $n=0$ there is nothing to prove. Next,
$$\begin{array}{l}
P_{k,n+1}(x,z)=P_{k,n}(x,z)\cdot (x-\beta_{k+n,1}y)\\[1em]\hspace{.5cm}
=\sum_a\Sha_{n-a,a}^{\uparrow k}\beta_{k,a}^{\phantom{\uparrow}}P_{0,a}^{\phantom{\uparrow}}(y,z)
P_{k,n-a}^{\phantom{\uparrow}}(x,y)^{\uparrow a}\cdot (x-\beta_{k+n,1}z)\\[1em]\hspace{.5cm}
=\sum_a\Sha_{n-a,a}^{\uparrow k}\beta_{k,a}^{\phantom{\uparrow}}P_{0,a}^{\phantom{\uparrow}}(y,z)
P_{k,n-a}^{\phantom{\uparrow}}(x,y)^{\uparrow a}\cdot \Bigl(
(x-\beta_{k+n-a,1}^{\uparrow a}y)+\beta_{k+n-a,1}^{\uparrow a}(y-\beta_{a,1}z)\Bigr)\\[1em]\hspace{.5cm}
=\sum_a\Sha_{n-a,a}^{\uparrow k}\Bigl(\beta_{k,a}^{\phantom{\uparrow}}P_{0,a}^{\phantom{\uparrow}}(y,z)
P_{k,n-a+1}^{\phantom{\uparrow}}(x,y)^{\uparrow a}+\beta_{k,a}^{\phantom{\uparrow}}\beta_{k+n-a,1}^{\uparrow a}
P_{0,a+1}^{\phantom{\uparrow}}(y,z)P_{k,n-a}(x,y)^{\uparrow a+1}\Bigr)\\[1em]\hspace{.5cm}
=\sum_a(\Sha_{n-a,a}^{\phantom{\uparrow}}+\Sha_{n-a+1,a-1}^{\phantom{\uparrow}}\beta_{n-a+1,1}^{\uparrow a-1})^{\uparrow k}
\beta_{k,a}^{\phantom{\uparrow}}P_{0,a}^{\phantom{\uparrow}}(y,z)P_{k,n-a+1}^{\phantom{\uparrow}}(x,y)^{\uparrow a}\\[1em]\hspace{.5cm}
=\sum_a\Sha_{n+1-a,a}^{\uparrow k}
\beta_{k,a}^{\phantom{\uparrow}}P_{0,a}^{\phantom{\uparrow}}(y,z)P_{k,n-a+1}^{\phantom{\uparrow}}(x,y)^{\uparrow a}\ .\end{array}$$
In the fourth equality we used (\ref{recpo}) to write $P_{k,n-a}(x,y)^{\uparrow a}(x-\beta_{k+n-a,1}^{\uparrow a}y)=P_{k,n+1-a}(x,y)^{\uparrow a}$;
then we moved $\beta_{k+n-a,1}^{\uparrow a}$ (the factor in front of $(y-\beta_{a,1}z)$) to the left through $P_{0,a}(y,z)P_{k,n-a}(x,y)^{\uparrow a}$
and it became $P_{0,a}(y,z)P_{k,n-a}(x,y)^{\uparrow a+1}$ in view of (\ref{br3}) and far commutativity; then, $(y-\beta_{a,1}z)$ commutes with
$P_{k,n-a}(x,y)^{\uparrow a+1}$ and $P_{0,a}(y,z)(y-\beta_{a,1}z)=P_{0,a+1}(y,z)$. In the fifth equality we separated the sums, used that
$\beta_{k,a}^{\phantom{\uparrow}}\beta_{k+n-a,1}^{\uparrow a}=\beta_{n-a,1}^{\uparrow k+a}\beta_{k,a+1}^{\phantom{\uparrow}}$ by (\ref{br5}), 
shifted the summation index in the second sum 
and rewrote the result as a single sum. In the sixth equality we used (\ref{shufrec}).\hfill $\Box$

\vskip .2cm
Since $P_{i,j}(x,0)=x^j$ and $P_{0,j}(0,z)=(-1)^j \omega_j z^j$, we find, evaluating  (\ref{robrbin}) at $y=0$, that
 \begin{equation} P_{k,n}(x,z)=  
 \sum_a (-1)^a\Sha_{n-a,a}^{\uparrow k}\beta_{k,a}^{\phantom{\uparrow}}\omega_a^{\phantom{\uparrow}} z^a x^{n-a}\ .\label{cobith}\end{equation}

\vskip .2cm
\noindent {\bf Braid Vandermonde identity.} Substituting (\ref{cobith}) into (\ref{recpo}) and collecting the powers of $x$, we obtain,
using (\ref{br4})-(\ref{br7}) and far commutativity, the braid version of the 
Vandermonde identity
\begin{equation}\Sha_{m+n-a,a}^{\uparrow k}\beta_{k,a}^{\phantom{\uparrow}}=\sum_{b,c:b+c=a} \Sha_{m-b,b}^{\uparrow k}
\Sha_{n-c,c}^{\uparrow m+k}\beta_{k,b}^{\phantom{\uparrow}}\beta_{m+k-b,c}^{\uparrow b}\ \ \ ,\ \ 
a,m,n,k\in {\mathbb Z}_{\geq 0} \ .\label{vandermonde1}\end{equation}
Since $\beta_{k,b}^{\phantom{\uparrow}}\beta_{m+k-b,c}^{\uparrow b}=\beta_{m-b,c}^{\uparrow b}\beta_{k,a}^{\phantom{\uparrow}}$
by (\ref{br5}),  eq.(\ref{vandermonde1}) can be derived from its particular case $k=0$,
\begin{equation}\Sha_{m+n-a,a}^{\phantom{\uparrow}}=\sum_{b,c:b+c=a} \Sha_{m-b,b}^{\phantom{\uparrow}}
\Sha_{n-c,c}^{\uparrow m}\beta_{m-b,c}^{\uparrow b}\ \ \ ,\ \ 
a,m,n\in {\mathbb Z}_{\geq 0} \ .\label{vandermonde1bis}\end{equation}
Eq.(\ref{vandermonde1bis}) generalizes the defining recursion for the shuffle elements: setting $n$ to 1, we reproduce (\ref{shufrec});
setting $m$ to 1, we reproduce (\ref{shufrec2}).

Another version of the Vandermonde identity is 
 \begin{equation}\omega_j^{\phantom{\uparrow}}\Sha_{e,c}^{\uparrow j}=\sum_a(-1)^a\Sha_{j-a,a}^{\phantom{\uparrow}}\omega_{j-a}^{\uparrow a}
\Sha_{e-a,c+j}^{\uparrow a}\ \beta_{a,c+j}^{\phantom{\uparrow}}\ \ \ ,\ \ 
e,c,j\in {\mathbb Z}_{\geq 0} \ .\label{vandermonde2}\end{equation}
It is proved by induction on $e+c$. For $e=0$ the relation (\ref{vandermonde2}) clearly holds. Thus, due to the initial conditions, it is 
enough to increase $e$, assuming that (\ref{vandermonde2}) holds for all smaller values of $e+c$. 
The combination ${\cal{X}}_{a,b;c}:=\Sha_{a,b}^{\uparrow c}\ \beta_{c,b}^{\phantom{\uparrow}}$ (entering the right hand side of
(\ref{vandermonde2})) verifies the recursion ${\cal{X}}_{a+1,b;c}={\cal{X}}_{a,b;c}+{\cal{X}}_{a+1,b-1;c}\beta_{a+c+1}^{\uparrow b-1}$ 
which straightforwardly implies the induction step.
   
\vskip .2cm
\noindent {\tt Notes.} 1.  The lower labels of many (but not all) elements carry information about the 
number of strands needed to define these elements (and therefore the number of copies of the vector
space in a local representation of the braid group tower). For example, $\Sha_{i,j},P_{i,j}(x,y),\beta_{i,j}\in {\mathfrak{k}}B_{i+j}$ and
$\omega_i\in {\mathfrak{k}}B_{i}$; this should not be confused with a different meaning of the lower label $j$ of $\sigma_j$.

\noindent 2. The case $k=0$ of the (\ref{robrbin}), (\ref{cobith}) and (\ref{vandermonde1}) was discussed in \cite{A}. 

\section{Braiding on the space of tensors}
\label{braidingtensors}

Extending the exchange rules (\ref{br2}) to all, not necessarily decomposable, tensors, we obtain
the standard braiding 
 \begin{equation}\x{[k]}\yu{[l]}{k}=\b_{k,l} \y{[l]}\xu{[k]}{l}\ ,\label{diag-braiding}\end{equation}
on the family $\{ V^{\otimes k}\}_{k=0}^\infty$ of vector spaces (see, e.g., \cite{K}). The 
corresponding braiding of the space $T(V)$ we call {\it diagonal}. The diagonal braiding of tensors is used, for example, in the construction
of the $q$-Minkowski vectors as bi-spinors \cite{OSWZ,OSWZ2} as well as in the $q$-versions of the 
accidental isomorphisms of semi-simple Lie groups \cite{JO}.

We shall construct another natural braiding $T(\hat{R})$ on the tensor space $T(V):=\bigoplus_k V^{\otimes k}$. "Natural" here means 
(i) functorial with respect to $V$; (ii) there is a quantum group associated to the 
braided space $V$ (see, e.g. \cite{FRT}). Its action naturally extends to $T(V)$. The braiding $T(\hat{R})$ is 
covariant with respect to this quantum group. 

Viewed differently, the above covariance reflects a particular property of our solution: the building blocks of the braiding $T(\hat{R})$ are certain 
polynomials in the original $R$-matrix; the Yang-Baxter system of equations can be formulated on the universal level of the braid group.
Our solution belongs to the ring of the monoid of positive braids. Questions related to the uniqueness of our solution will be 
discussed in \cite{GO}. A version of the braiding on the tensor space in which
the building blocks of  $T(\hat{R})$ are filled modulo a given natural number $N$ will be considered in \cite{GO2}.

Due to the canonical isomorphisms $V^{\otimes l}\otimes V^{\otimes k}\rightarrow V^{\otimes (l+k)}$, we can consider
a more general Ansatz for a braiding on $T(V)$, for which the right hand side of (\ref{diag-braiding}) acquires other terms
$\y{[l']}\xu{[k']}{l'}$ with $k'+l'=k+l$. Our aim is to study such braidings.

Denote by $\A{a}{b}{c}$ the submatrix appearing in front of a term $\y{[k]}\xu{[b+c-k]}{k}$.
Thus the braiding $T(\hat{R})$, understood as the exchange rule, reads
\begin{equation}\label{braiding_mat1}
\x{[b]}\yu{[c]}{b}=\sum_{k=0}^{b+c}\A{b}{c}{k}~\y{[k]}\xu{[b+c-k]}{k}.
\end{equation}
The full form of $\A{b}{c}{k}~\y{[k]}\xu{[b+c-k]}{k}$ is 
${\A{b}{c}{k}}^{i_1\dots i_{b+c}}_{j_1\dots j_{b+c}}y^{j_1\dots j_k}x^{j_{k+1}\dots j_{b+c}}$, 
the summation in repeated indices is assumed.

As for the standard (Drinfeld--Jimbo, \cite{D,J}) constant $R$-matrices, it turns out that a braiding 
with all $\A{b}{c}{k}$ non-vanishing does not exist. We require the $R$-matrix for $T(V)$ to 
be block-triangular (like in the construction of the orthogonal and symplectic $R$-matrices
in \cite{O}). In other words we restrict ourselves to braidings with $\A{b}{c}{k}=0$ if 
$k>c$. Then the summation in (\ref{braiding_mat1}) shortens to
\begin{equation}\label{relcomtens'}\x{[b]}\yu{[c]}{b}=
\sum_{k=0}^{c}\A{b}{c}{c-k}~\y{[c-k]}\xu{[b+k]}{c-k}\ .\end{equation}
For example, 
$$\x{[3]}\yu{[2]}{3}=\A{3}{2}{2}~\y{[2]}\xu{[3]}{2}+
\A{3}{2}{1}~\y{[1]}\xu{[4]}{1}+\A{3}{2}{0}~\y{[0]}\xu{[5]}{0}\ .$$
The operation $\ ^{\uparrow 0}$ is the identity and can be omitted, $\xu{[5]}{0}=\x{[5]}$.

The submatrix $\A{1}{1}{1}$ is an endomorphism of $V\otimes V$. With our triangular Ansatz, this endomorphism must verify
the Yang-Baxter equation. We identify $\A{1}{1}{1}$ with the initial braiding on $V$. We shall see that our solution contains
two "parameters", the Yang-Baxter operator $\hat{R}=\A{1}{1}{1}$ and the scalar $q=\A{0}{0}{0}$ (by convention, 
$\x{[0]}\in V^{\otimes 0}\simeq \k$).

To be a braiding, $T(\hat{R})$ must satisfy the Yang-Baxter equation. As for vectors, we shall understand the Yang--Baxter 
equation for an operator $\hat{R}$ as the equality of two different reorderings of $x^\bullet y^\bullet z^\bullet$ (using 
the rule (\ref{relcomtens'}) for $x^\bullet y^\bullet$, $x^\bullet z^\bullet$ and $y^\bullet z^\bullet$)
to the form $z^\bullet y^\bullet x^\bullet$. 

We have
\begin{equation}
\x{[a]}(\yu{[b]}{a}\zu{[c]}{a+b})=\sum_{d,e,f}\Au{b}{c}{d}{a}\A{a}{d}{e}\Au{a+d-e}{b+c-d}{f}{e}\z{[e]} \yu{[f]}{e} \xu{[a+b+c-e-f]}{e+f}
\end{equation}
and
\begin{equation}
(\x{[a]}\yu{[b]}{a})\zu{[c]}{a+b}=\sum_{u,v,w}\A{a}{b}{u}\Au{a+b-u}{c}{v}{u}\A{u}{v}{w}\z{[w]} \yu{[u+v-w]}{w} \xu{[a+b+c-u-v]}{u+v}.
\end{equation}
Equating terms, we arrive at a system
\begin{equation}\sum_{i}\Au{b}{c}{i}{a}\A{a}{i}{e}\Au{a+i-e}{b+c-i}{f}{e}=
\sum_{j}\A{a}{b}{j}\Au{a+b-j}{c}{f+e-j}{j}\A{j}{f+e-j}{e}\ .\label{sytso}\end{equation}
Here the summation is over the terms which are non zero. This implies that $a,b,c,e,f\in {\mathbb Z}_{\geq 0}$; moreover, in the left hand 
side the summation is over $i$ such that $\max(0,e)\leq i\leq \min (c,b+c-f)$ and, in the right hand side, the summation is over $j$ such that
$\max(0,f+e-c)\leq j\leq \min (f,b)$. 

The system (\ref{sytso}) can be considered on the "universal" level, as a system of equations in the braid group ring. We shall take this
point of view and construct a universal solution of the system (\ref{sytso}), with $\A{a}{b}{c}\in \mathfrak{k} B_{a+b}$.  
For a given braiding $\hat{R}$ on the vector space $V$, the braiding $T(\hat{R})$ on $T(V)$ is obtained by taking the local representation 
$\rho_{q\hat{R}}$ of the braid group tower ($\rho_{q\hat{R}}$ sends $\sigma_i$ to $q\hat{R}_i$).  

\vskip .2cm
\noindent {\bf Theorem}. Let $p_{a,b}:=P_{a,b}(1,q^{-2})$. The following elements
\begin{equation}
\A{a}{b}{c}=q^{1-a-c}\ \Sha_{b-c,c}^{\uparrow a}\beta_{a,c}^{\phantom{\uparrow}}p_{a,b-c}^{\uparrow c}\ ,
\end{equation}
provide a universal solution of the system (\ref{sytso}). 

\vskip .2cm
\noindent {\it Proof}. The left hand side of (\ref{sytso}) multiplied by $q^{2a+b+f-3}$ reads:
$$\sum_i q^{-2i}\Sha_{c-i,i}^{\uparrow a+b}\beta_{b,i}^{\uparrow a}p_{b,c-i}^{\uparrow a+i}
\Sha_{i-e,e}^{\uparrow a}\beta_{a,e}^{\phantom{\uparrow}}p_{a,i-e}^{\uparrow e}
\Sha_{b+c-f-i,f}^{\uparrow a+i}\beta_{a+i-e,f}^{\uparrow e}p_{a+i-e,b+c-f-i}^{\uparrow e+f}\ .$$
The term $p_{a,i-e}^{\uparrow e}$ commutes with $\Sha_{b+c-f-i,f}^{\uparrow a+i}$ (far commutativity). When $p_{a,i-e}^{\uparrow e}$ 
moves to the right through $\beta_{a+i-e,f}^{\uparrow e}$ it becomes $p_{a,i-e}^{\uparrow e+f}$ by (\ref{br3}). Next, by (\ref{recpo}), 
$p_{a,i-e}^{\uparrow e+f}p_{a+i-e,b+c-f-i}^{\uparrow e+f}=p_{a,b+c-f-e}^{\uparrow e+f}$ and it does not depend on the summation
index. After, the term $\Sha_{i-e,e}^{\uparrow a}$ commutes with $p_{b,c-i}^{\uparrow a+i}$ (far commutativity); then 
$\Sha_{i-e,e}^{\uparrow a}$ moves to the left through $\beta_{b,i}^{\uparrow a}$ and becomes 
$\Sha_{i-e,e}^{\uparrow a+b}$ by (\ref{br3}); by (\ref{shaiden}), we have 
$\Sha_{c-i,i}^{\uparrow a+b}\Sha_{i-e,e}^{\uparrow a+b}=\Sha_{c-e,e}^{\uparrow a+b}\Sha_{c-i,i-e}^{\uparrow a+b+e}$ and 
the left hand side of (\ref{sytso}) becomes 
\begin{equation}\Sha_{c-e,e}^{\uparrow a+b}\Bigl(\sum_i q^{-2i}\Sha_{c-i,i-e}^{\uparrow a+b+e}\beta_{b,i}^{\uparrow a}p_{b,c-i}^{\uparrow a+i}
\beta_{a,e}^{\phantom{\uparrow}}\Sha_{b+c-f-i,f}^{\uparrow a+i}\beta_{a+i-e,f}^{\uparrow e}\Bigr)
p_{a,b+c-f-e}^{\uparrow e+f}\ .\label{lhssysto1}\end{equation} 
The term $\beta_{a,e}$ commutes with $\Sha_{b+c-f-i,f}^{\uparrow a+i}$ (far commutativity). By (\ref{br5}), $\beta_{a,e}^{\phantom{\uparrow}}
\beta_{a+i-e,f}^{\uparrow e}=\beta_{i-e,f}^{\uparrow e+a}\beta_{a,e+f}$; the factor $\beta_{a,e+f}$ does not depend on 
the summation index and goes out from the sum. By (\ref{br4}), $\beta_{b,i}^{\uparrow a}=\beta_{b,e}^{\uparrow a}\beta_{b,i-e}^{\uparrow a+e}$. 
The factor $\beta_{b,e}^{\uparrow a}$ does not depend on the summation index, it commutes with $\Sha_{c-i,i-e}^{\uparrow a+b+e}$ 
 (far commutativity) and moves out from the sum to the left. The left hand side of (\ref{sytso}) takes now the form
\begin{equation}\Sha_{c-e,e}^{\uparrow a+b}\beta_{b,e}^{\uparrow a}\Bigl(\sum_i q^{-2i}\Sha_{c-i,i-e}^{\uparrow b+e}
\beta_{b,i-e}^{\uparrow e}p_{b,c-i}^{\uparrow i}\Sha_{b+c-f-i,f}^{\uparrow i}\beta_{i-e,f}^{\uparrow e}\Bigr)^{\uparrow a}
\beta_{a,e+f}^{\phantom{\uparrow}}p_{a,b+c-f-e}^{\uparrow e+f}\ .\label{lhssysto2}\end{equation}  

The right hand side of (\ref{sytso}) multiplied by $q^{2a+b+f-3}$ reads:
$$q^{-2e}\sum_j \Sha_{b-j,j}^{\uparrow a}\beta_{a,j}p_{a,b-j}^{\uparrow j}
\Sha_{c-e-f+j,e+f-j}^{\uparrow a+b}\beta_{a+b-j,f+e-j}^{\uparrow j}p_{a+b-j,c-e-f+j}^{\uparrow e+f}
\Sha_{f-j,e}^{\uparrow j}\beta_{j,e}p_{j,f-j}^{\uparrow e}\ .$$
The term $p_{a,b-j}^{\uparrow j}$ commutes with $\Sha_{c-e-f+j,e+f-j}^{\uparrow a+b}$  (far commutativity); when $p_{a,b-j}^{\uparrow j}$ 
moves to the right through 
$\beta_{a+b-j,f+e-j}^{\uparrow j}$, it becomes $p_{a,b-j}^{\uparrow e+f}$ by (\ref{br3}); then, by (\ref{recpo}) we have 
$p_{a,b-j}^{\uparrow e+f}p_{a+b-j,c-e-f+j}^{\uparrow e+f}=p_{a,b+c-e-f}^{\uparrow e+f}$; the term $p_{a,b+c-e-f}^{\uparrow e+f}$
moves to the right out of the sum without changes (far commutativity). Next, the term $\Sha_{f-j,e}^{\uparrow j}$ moves to the left 
through $\beta_{a+b-j,f+e-j}^{\uparrow j}$, becoming $\Sha_{f-j,e}^{\uparrow a+b}$ by (\ref{br3}); 
we have $\Sha_{c-e-f+j,e+f-j}^{\uparrow a+b}\Sha_{f-j,e}^{\uparrow a+b}=
\Sha_{c-e,e}^{\uparrow a+b}\Sha_{c-e-f+j,f-j}^{\uparrow a+b+e}$  by (\ref{shaiden}); the term 
$\Sha_{c-e,e}^{\uparrow a+b}$ moves now without changes to
the very left (far commutativity) and the right hand side becomes 
\begin{equation}\Sha_{c-e,e}^{\uparrow a+b}\Bigl(q^{-2e}\sum_j\Sha_{b-j,j}^{\uparrow a}\beta_{a,j}^{\phantom{\uparrow}}
\Sha_{c-e-f+j,f-j}^{\uparrow a+b+e}\beta_{a+b-j,f+e-j}^{\uparrow j}\beta_{j,e}^{\phantom{\uparrow}}p_{j,f-j}^{\uparrow e}\Bigr)
p_{a,b+c-f-e}^{\uparrow e+f}\ .\label{rhssysto1}\end{equation}  
The term $\beta_{a,j}^{\phantom{\uparrow}}$ commutes with $\Sha_{c-e-f+j,f-j}^{\uparrow a+b+e}$ (far commutativity); then, 
by (\ref{br5}), $\beta_{a,j}^{\phantom{\uparrow}}\beta_{a+b-j,f+e-j}^{\uparrow j}=
\beta_{b-j,e+f-j}^{\uparrow a+j}\beta_{a,e+f}^{\phantom{\uparrow}}$; the term $\beta_{a,e+f}^{\phantom{\uparrow}}$ does not depend
on the summation index; we move it out of the sum to the right; when it moves through $\beta_{j,e}^{\phantom{\uparrow}}p_{j,f-j}^{\uparrow e}$,
this expression transforms into $\beta_{j,e}^{\uparrow a}p_{j,f-j}^{\uparrow a+e}$. We then rewrite: $\beta_{b-j,e+f-j}^{\uparrow a+j}
\beta_{j,e}^{\uparrow a}=\beta_{b,e}^{\uparrow a}\beta_{b-j,f-j}^{\uparrow a+j+e}$ by (\ref{br5}) and move the term 
$\beta_{b,e}^{\uparrow a}$ out of the sum to the left; it commutes with $\Sha_{c-e-f+j,f-j}^{\uparrow a+b+e}$  (far commutativity) and 
transforms $\Sha_{b-j,j}^{\uparrow a}$ into $\Sha_{b-j,j}^{\uparrow a+e}$ by (\ref{br3}). The right hand side of (\ref{sytso}) 
takes the form 
\begin{equation}\Sha_{c-e,e}^{\uparrow a+b}\beta_{b,e}^{\uparrow a}\Bigl(q^{-2e}\sum_j\Sha_{b-j,j}^{\uparrow e}
\Sha_{c-e-f+j,f-j}^{\uparrow b+e}\beta_{b-j,f-j}^{\uparrow j+e}p_{j,f-j}^{\uparrow e}\Bigr)^{\uparrow a}
\beta_{a,e+f}^{\phantom{\uparrow}}p_{a,b+c-f-e}^{\uparrow e+f}\ .\label{rhssysto2}\end{equation}    
 
Comparing (\ref{lhssysto2}) with (\ref{rhssysto2}), we see that the theorem will follow from the equality of 
\begin{equation}
\sum_i q^{-2i}\Sha_{c-i,i-e}^{\uparrow b}
\beta_{b,i-e}^{\phantom{\uparrow}}p_{b,c-i}^{\uparrow i-e}\Sha_{b+c-f-i,f}^{\uparrow i-e}\beta_{i-e,f}^{\phantom{\uparrow}}=
q^{-2e}\sum_j \Sha_{b-j,j}^{\phantom{\uparrow}}\Sha_{c-e-f+j,f-j}^{\uparrow b}\beta_{b-j,f-j}^{\uparrow j}p_{j,f-j}^{\phantom{\uparrow}}
\ .\label{sytso2}\end{equation}
Substitute (\ref{cobith}) in the left hand side of (\ref{sytso2}):
\begin{equation}
\sum_{i,t} (-1)^t q^{-2i-2t}\Sha_{c-i,i-e}^{\uparrow b}
\beta_{b,i-e}^{\phantom{\uparrow}}\Sha_{c-i-t,t}^{\uparrow b+i-e}\beta_{b,t}^{\uparrow i-e}\omega_t^{\uparrow i-e}
\Sha_{b+c-f-i,f}^{\uparrow i-e}\beta_{i-e,f}^{\phantom{\uparrow}}
\ .\label{sytso2l1}\end{equation}
Move $\beta_{b,i-e}^{\phantom{\uparrow}}$ to the right through $\Sha_{c-i-t,t}^{\uparrow b+i-e}$ (far commutativity); use (\ref{br4}) to 
write $\beta_{b,i-e}^{\phantom{\uparrow}}\beta_{b,t}^{\uparrow i-e}=\beta_{b,i+t-e}^{\phantom{\uparrow}}$ and (\ref{shaiden}) to 
write $\Sha_{c-i,i-e}^{\uparrow b}\Sha_{c-i-t,t}^{\uparrow b+i-e}=\Sha_{c-i-t,i+t-e}^{\uparrow b}\Sha_{t,i-e}^{\uparrow b}$; the term
$\Sha_{t,i-e}^{\uparrow b}$ moves now to the right through $\beta_{b,i+t-e}^{\phantom{\uparrow}}$ becoming 
$\Sha_{t,i-e}^{\phantom{\uparrow}}$. Replace the summation index $t$ by $r=i+t$:
\begin{equation}
\sum_{i,r} (-1)^{r+e} q^{-2r}\Sha_{c-r,r-e}^{\uparrow b}
\beta_{b,r-e}^{\phantom{\uparrow}}\Bigl\{ (-1)^{i-e}\Sha_{r-i,i-e}^{\phantom{\uparrow}}\omega_{r-i}^{\uparrow i-e}
\Sha_{b+c-f-i,f}^{\uparrow i-e}\beta_{i-e,f}^{\phantom{\uparrow}}\Bigr\}
\ .\label{sytso2l2}\end{equation}
The sum over $i$ (terms in braces) is taken with the help of (\ref{vandermonde2}) and the final form of the left hand side of (\ref{sytso2}) reads
\begin{equation}
\sum_{r} (-1)^{r+e} q^{-2r}\Sha_{c-r,r-e}^{\uparrow b}
\beta_{b,r-e}^{\phantom{\uparrow}}\omega_{r-e}^{\phantom{\uparrow}}
\Sha_{b+c-e-f,e+f-r}^{\uparrow r-e}
\ .\label{sytso2l3}\end{equation}
Substitute (\ref{cobith}) in the right hand side of (\ref{sytso2}):
\begin{equation}
q^{-2e}\sum_{j,s}(-1)^s q^{-2s} \Sha_{b-j,j}^{\phantom{\uparrow}}\Sha_{c-e-f+j,f-j}^{\uparrow b}\beta_{b-j,f-j}^{\uparrow j}
\Sha_{f-j-s,s}^{\uparrow j}\beta_{j,s}^{\phantom{\uparrow}}\omega_s^{\phantom{\uparrow}}
\ .\label{sytso2r1}\end{equation}
The element $\Sha_{f-j-s,s}^{\uparrow j}$ moves to the left through $\beta_{b-j,f-j}^{\uparrow j}$ becoming
$\Sha_{f-j-s,s}^{\uparrow b}$ by (\ref{br3}). Then we write $\Sha_{c-e-f+j,f-j}^{\uparrow b}\Sha_{f-j-s,s}^{\uparrow b}=
\Sha_{c-e-s,s}^{\uparrow b}\Sha_{c-e-f+j,f-j-s}^{\uparrow b+s}$ by (\ref{shaiden}) and $\beta_{b-j,f-j}^{\uparrow j}\beta_{j,s}^{\phantom{\uparrow}}
=\beta_{b,s}^{\phantom{\uparrow}}\beta_{b-j,f-j-s}^{\uparrow j+s}$ by (\ref{br5}). The shuffle element $\Sha_{c-e-s,s}^{\uparrow b}$
moves through $\Sha_{b-j,j}^{\phantom{\uparrow}}$ (far commutativity) to the left; the term $\beta_{b,s}^{\phantom{\uparrow}}$ moves to the left
through the product of two shuffles, the shuffle $\Sha_{c-e-f+j,f-j-s}^{\uparrow b+s}$ does not change (far commutativity) while the
shuffle $\Sha_{b-j,j}^{\phantom{\uparrow}}$ becomes $\Sha_{b-j,j}^{\uparrow s}$. The right hand side of (\ref{sytso2}) reads now
\begin{equation}
q^{-2e}\sum_{j,s}(-1)^s q^{-2s} \Sha_{c-e-s,s}^{\uparrow b}\beta_{b,s}^{\phantom{\uparrow}}
\Bigl\{ \Sha_{b-j,j}^{\uparrow s}\Sha_{c-e-f+j,f-j-s}^{\uparrow b+s}\beta_{b-j,f-j-s}^{\uparrow j+s}\Bigr\}
\omega_s^{\phantom{\uparrow}}
\ .\label{sytso2r2}\end{equation}
The sum in $j$ (terms in braces) is taken with the help of (\ref{vandermonde1}), it equals $\Sha_{b+c-e-f,f-s}^{\uparrow s}$, which
far commutes with $\omega_s^{\phantom{\uparrow}}$ and we arrive at the same expression (\ref{sytso2l3}). The proof is completed.\hfill $\Box$

\end{document}